\documentclass[11pt]{article}

\usepackage{amsmath,amssymb}
\usepackage[latin1]{inputenc}
\usepackage[dvips]{graphicx}

\topskip=\baselineskip  \headsep=0in \headheight=0in
\usepackage{amsmath,amssymb}

\newtheorem{theorem}{Theorem}[section]
\newtheorem{lemma}[theorem]{Lemma}
\newtheorem{proposition}[theorem]{Proposition}

\newtheorem{corollary}[theorem]{Corollary}

\newcommand{\qed}{\enspace\vrule  height6pt  width4pt  depth2pt}

\def\thetheorem{\arabic{theorem}}

\newenvironment{proof}{\par\noindent{\bf Proof.}}{$\qed$\par\bigskip}

\newcommand{\Aut}{\mbox{\rm Aut}}

\newcommand{\Imagen}{\mbox{\rm Im }}

\newcommand{\End}{\mbox{\rm End}}
\newcommand{\Hom}{\mbox{\rm Hom}}

\newcommand{\inv}{^{-1}}

\newcommand{\Gal}{{\rm Gal}}

\newcommand{\GL}{{\rm GL}}

\newcommand{\tr}{{\rm tr}}
\newcommand{\Ind}{{\rm Ind}}
\newcommand{\matriz}[1]{\begin{array} #1 \end{array}}

\newcommand{\GEN}[1]{\langle #1 \rangle}

\title{Every projective Schur algebra is Brauer equivalent to a radical abelian algebra}
\author{Eli Aljadeff\thanks{This research was carried
while the first author was visiting the University of Murcia during the spring semester of 2005/6. He thanks the
Department of Mathematics for the hospitality and Fundaci\'{o}n S\'{e}neca of Murcia for their financial support.}
\hspace{0.05cm} and \'{A}ngel del R\'{\i}o\thanks{The second author has been partially supported by D.G.I. of Spain and
Fundaci\'{o}n S\'{e}neca of Murcia.}}

\begin{document}

\maketitle

\begin{abstract}
We prove that any projective Schur algebra over a field $K$ is equivalent in $Br(K)$ to a radical abelian algebra. This
was conjectured in 1995 by Sonn and the first author of this paper. As a consequence we obtain a characterization of
the projective Schur group by means of Galois cohomology. The conjecture was known for algebras over fields of positive
characteristic. In characteristic zero the conjecture was known for algebras over fields with an Henselian valuation
over a local or global field of characteristic zero.
\end{abstract}

\section{Introduction}

One of the main theorems on the Schur subgroup of the Brauer group
is the Brauer-Witt Theorem.

\medskip
\noindent{\bf Theorem} (Brauer-Witt). Every Schur algebra is
Brauer equivalent to a cyclotomic algebra. \vspace{0.1cm}

Recall that a (finite dimensional) central simple $K$-algebra $A$
is called a Schur algebra if it is spanned over $K$ by a finite
group of unit elements. Equivalently, $A$ is a Schur algebra if it
is the homomorphic image of a group algebra $KG$ for some finite
group $G$. The subgroup of $Br(K)$ generated by (in fact
consisting of) classes that are represented by Schur algebras is
called the Schur group of $K$ and is denoted by $S(K)$ \cite{Y}.
There is a natural way to construct Schur algebras, namely the
cyclotomic algebras. A cyclotomic algebra over $K$ is a crossed
product $A=(L/K,\alpha)$ where $L=K(\zeta)/K$ is a cyclotomic
field extension and $\alpha \in H^{2}(\Gal(L/K),L^{\times})$ is
represented by a $2$-cocycle $t$ whose values are roots of unity.
Clearly, if $A =\bigoplus_{\sigma\in \Gal(L/K)} Lu_{\sigma}$ where
$u_{\sigma}u_{\tau}=t(\sigma,\tau)u_{\sigma\tau}$ then the group
of units generated by $\zeta$ and the $u_{\sigma}$'s is finite and
spans $A$ as a vector space over $K$. Hence $A$ is a Schur
algebra. The Brauer-Witt Theorem states that this natural
construction amounts for all the Schur algebras up to Brauer
equivalence.

Using the Brauer-Witt theorem, one can characterize the Schur
group by means of Galois cohomology. Let $K_{s}$ be a separable
closure of $K$ and let $K_{cyc}$ be the maximal cyclotomic
extension of $K$ (contained in $K_{s}$). Denote by
$G_{s}=\Gal(K_{s}/K)$ and $G_{cyc}=\Gal(K_{cyc}/K)$ the
corresponding Galois groups.

\medskip
\noindent{\bf Corollary}. Let $\mu$ be the group of roots of unity in $K_{cyc}$ viewed as a $G_{cyc}$ module. Consider
the following commutative diagram where $i_{*}: H^{2}(G_{cyc},\mu)\rightarrow H^{2}(G_{cyc}, K_{cyc}^{\times})$ is the map
induced by the inclusion $i:\mu \hookrightarrow K_{cyc}^{\times}$ and $\phi $ is the map induced by the cyclotomic algebra
construction.
    $$\matriz{{rccccc}
    H^{2}(G_{cyc}, \mu) \stackrel{i_*}{\twoheadrightarrow}  & \Imagen
i_*            & \hookrightarrow &
                                                H^{2}(G_{cyc},
K_{cyc}^{\times}) & \rightarrow & H^{2}(G_{K}, K_{s}^{\times}) \\
                        \phi \searrow
      & \bar{\phi} \downarrow\wr  &                 &
                                                \downarrow\wr
     &                   & \downarrow\wr \\
                                                      &  S(K)
        & \hookrightarrow &
                                                Br(K_{cyc}/K)
     & \hookrightarrow   & Br(K)}$$
Then $\bar{\phi}$ is an isomorphism. \medskip

In 1978 projective analogues of these constructions were introduced by Lorenz and Opolka \cite{LO}. Recall that a
central simple $K$-algebra $A$ is called projective Schur if it is spanned by a group of units which is finite modulo
its centre. Equivalently, $A$ is projective Schur if it is the homomorphic image of a twisted group algebra $K^{t}G$
for some finite group $G$ and $2$-cocycle $t$. As in the Schur case, the classes of $Br(K)$ represented by projective
Schur algebras form a group. This is the projective Schur group of $K$ and is denoted by $PS(K)$. One of the main
points for introducing these notions is that every symbol algebra is projective Schur in an obvious way. This
observation already shows that the projective Schur group is $large$ if roots of units are present \cite{NO}. Indeed
the Merkurjev-Suslin theorem implies that if $K$ contains all roots of unity then $PS(K)=Br(K)$. Interestingly, as
observed by Lorenz and Opolka, also for fields which have only few roots of unity, the projective Schur group may be
$large$. They showed that if $K$ is a local or global field then $PS(K)=Br(K)$ \cite{LO}. Based on these examples it
was conjectured that $PS(K)=Br(K)$ for arbitrary fields but this turned out to be false as shown in \cite{AS2} (e.g.
$K=k(x)$, $k$ any number field).

The main objective was then to find a description of $PS(K)$ in terms of Galois cohomology. As in the Schur case, also
here there is a natural way to construct projective Schur algebras. These are the radical or radical abelian algebras.
Before we recall their definition let us introduce the following notation: Given any field extension $L/K$ we denote by
$Rad_K(L)$ the group of all units in $L$ which are of finite order modulo $K^{\times}$. A field extension $L/K$ is said to be
radical if $L=K(Rad_K(L))$\footnote{Note that this definition is more restrictive than the usual definition of radical
extension in Galois Theory}. Now, a central simple $K$-algebra is radical if it is a crossed product $A=(L/K,\alpha)$
where $L/K$ is a Galois radical extension and $\alpha \in Rad_K(L)$ (by this we mean that the cohomology class $\alpha$
may be represented by a $2$-cocycle $t$ with $t(g,h)\in Rad_K(L)$ for every $g,h\in \Gal(L/K))$. If the extension $L/K$
is abelian, we say that $A$ is radical abelian. Clearly, if $\{u_{g}:g\in Gal(L/K)\}$ is an $L$-basis with
$u_gu_h=t(g,h)u_{gh}$ ($t$ as above) then the multiplicative group $\GEN{Rad_K(L),u_g: g\in \Gal(L/K)}$ spans $A$ as a
vector space over $K$ and is finite modulo its centre. Hence, every radical algebra is projective Schur.

The "Brauer-Witt Conjecture" for projective Schur algebras says that every projective Schur algebra is Brauer
equivalent to a radical (abelian) algebra (see \cite{AS1}). In \cite{AS-E} the conjecture was proved for fields of
positive characteristic. In characteristic zero the conjecture was established for Henselian valued fields whose
residue field is a local or global field of characteristics zero (e.g. iterated Laurent series over a number field)
(see \cite{ASW}).

As the title suggests, the goal of this paper is to prove the Brauer-Witt conjecture for projective Schur algebras over
arbitrary fields of characteristic 0.

\begin{theorem}\label{PBW}
Every projective Schur algebra is Brauer equivalent to a radical
abelian algebra.
\end{theorem}

As a consequence we obtain a cohomological interpretation of $PS(K)$. Let $K_{rad,ab}$ be the maximal radical abelian
extension of $K$ (contained in $K_{s}$) and let $G_{rad,ab}=\Gal(K_{rad,ab}/K)$ be the Galois group.

\begin{corollary}
Consider the following commutative diagram where $i_{*}: H^{2}(G_{rad,ab},Rad_K(K_{rad,ab}))\rightarrow
H^{2}(G_{rad,ab}, K_{rad,ab}^{\times})$ is the map induced by the inclusion $i:Rad_K(K_{rad,ab}) \hookrightarrow K_{rad,ab}^{\times}$
and $\phi$ is the map induced by the radical algebra construction.
    $$\matriz{{rccccc}
    H^{2}(G_{rad,ab}, Rad_K(K_{rad,ab})) \stackrel{i_*}{\twoheadrightarrow}  & \Imagen i_*
            & \hookrightarrow & H^{2}(G_{rad,ab}, K_{rad,ab}^{\times}) & \rightarrow & H^{2}(G_{K}, K_{rad,ab}^{\times}) \\
                         \phi \searrow
      & \bar{\phi} \downarrow\wr  &                 &
                                                \downarrow\wr
     &                   & \downarrow\wr \\
                                                      &  PS(K)
         & \hookrightarrow &
                                                Br(K_{rad,ab}/K)
     & \hookrightarrow   & Br(K)}$$
Then $\bar{\phi}$ is an isomorphism.
\end{corollary}

The proof of Theorem~\ref{PBW} uses the idea of the proof of the Brauer-Witt Theorem as it appears in \cite{Y}. The
main difficulty here comes from a fundamental difference between cyclotomic and radical algebras. While being a root of
unity is a notion which is independent of the field considered, being radical depends on the ground field. This
obstruction appears in the use of the corestriction map.
To overcome this difficulty it is convenient to generalize the notions of projective Schur and radical algebra.

Let $K/k$ be a field extension and $A$ a central simple $K$-algebra. We say that $A$ is projective Schur over $k$ if it
is an epimorphic image of a twisted group algebra $k^t G$ of a finite group $G$. Equivalently $A$ is projective Schur
over $k$ if it is generated over $k$ by a group of units $\Gamma$ which is finite modulo $k^{\times}$ (that is
$[\Gamma:\Gamma\cap k^{\times}]<\infty$). We denote such an algebra by $k(\Gamma)$. We say that $A$ is radical over $k$ if $A$
is a crossed product $(L/K,\alpha)$ where $L/k$ is a (non necessarily Galois) radical extension and $\alpha \in
Rad_k(L)$.


Clearly, an algebra $A$ which is radical over $k$ is projective
Schur over $k$ in a natural way. In the theorem below we show that
the generating group $\Gamma$ that appears in this way is
supersolvable. In fact the class of supersolvable groups is the
``precise'' class of groups needed to represent radical algebras.


\begin{theorem}\label{RadSS}
Let $K/k$ be a field extension and $A$ a central simple
$K$-algebra.
\begin{enumerate}
\item[(a)] If $A$ is radical over $k$ then $A$ is generated by a
supersolvable group of units which is finite modulo $k$.
\item[(b)] If $A$ is generated by a supersolvable group of units
which is finite modulo $k^\times$, then $A$ is isomorphic to a full matrix algebra of an algebra which is radical
abelian over $k$.
\end{enumerate}
\end{theorem}

Theorem~\ref{RadSS} is proved in Section 2. Note that Theorem~\ref{RadSS}(b) implies Theorem~\ref{PBW} for projective
Schur algebras $k(\Gamma)$ with $\Gamma$ supersolvable. The main task is then to show that every projective Schur
algebra is equivalent to a projective Schur algebra which is generated by a supersolvable group. For fields of
characteristic zero this is contained in the next theorem. It implies Theorem~\ref{PBW} for fields of characteristic
zero. (For fields of positive characteristic, Theorem~\ref{PBW} was proved in \cite{AS-E}.)

\begin{theorem}\label{Main}
Let $K/k$ be a extension of fields of characteristic $0$ and $A$ a
central simple $K$-algebra. Then the following conditions are
equivalent.
\begin{enumerate}
\item
$A$ is equivalent (in $Br(K)$) to a central simple $K$-algebra which is projective Schur over $k$.
\item
$A$ is equivalent to a central simple $K$-algebra generated by a supersolvable group of units which is finite modulo
$k^\times$ (i.e. equivalent to an epimorphic image of a twisted group algebra $k^t G$ for $G$ a finite supersolvable
group).
\item
$A$ is equivalent to a central simple $K$-algebra which is radical over $k^{\times}$.
\item
$A$ is equivalent to a central simple $K$-algebra which is radical abelian over $k$.
\end{enumerate}
\end{theorem}

Implications (1 $\Leftarrow$ 2 $\Leftrightarrow$ 3 $\Leftrightarrow$ 4) are either obvious or follow from
Theorem~\ref{RadSS} (note that Brauer equivalence is basically not needed for these implications). To complete the
proof of Theorem~\ref{Main} one needs to show (1 $\Rightarrow$ 3). An important part of this proof is in
Proposition~\ref{RedE} which uses character theory of projective representation.


\section{Supersolvable groups}\label{SectSS}
\def\thetheorem{\thesection.\arabic{theorem}}

In this section we prove Theorem~\ref{RadSS}.

\medskip

\noindent{{\bf Proof of Theorem~\ref{RadSS}(a)}. Let $A=(L/K,\alpha=[t])=\oplus_{\sigma\in G} Lu_{\sigma}$ be radical
over $k$ with $G=\Gal(L/K)$. Let $L=k(\alpha_1,\dots,\alpha_k)$ with each $\alpha_i$ finite modulo $k^{\times}$. One may
assume that $t(\sigma,\tau) \in \GEN{\alpha_1,\dots,\alpha_n}$ for each $\sigma,\tau\in G$.

Let $E_0$ be a maximal cyclotomic extension of $K$ contained in
$L$. For each $i=1,\dots,k$, let $n_i$ be the order of $\alpha_i$
modulo $E_0^{\times}$ and $m_i$ the degree of the minimal
polynomial $p_i$ of $\alpha_i$ over $E_0$. Every conjugate of
$\alpha_i$ over $E_0$ is of the form $\alpha_i\zeta$ for some root
of unity $\zeta$ in $L$ and thus the coefficient of degree $0$ of
$p_i$ is $\alpha_i^{m_i} \gamma$ for some root of unity $\gamma\in
L$. By the maximality of $E_0$, $\gamma$ and each $\zeta$ belong
to $E_0$. Thus $\alpha_i^{m_i}\in E_0$ and therefore $m_i\ge n_i$.
This implies that $p_i=X^{n_i}-\alpha_i^{n_i}$ and so $E_0$
contains a primitive $n_i$-th root of unity. This shows that
$L/E_0$ is a Kummer extension and that $E_0$ contains a primitive
$n$-th root of unity $\xi$, where $n$ is the least common multiple
of the $n_i$'s.

Let $\Gamma$ be the subgroup of the group of units of $A$ generated by $\xi$, the $\alpha_i$'s and the $u_{\sigma}$'s.
Clearly $A$ is generated by $\Gamma$ over $k$ and $[\Gamma:\Gamma\cap k^{\times}]<\infty$. We show that $\Gamma$ is
supersolvable. Clearly $C_1=\GEN{\xi}$ is a cyclic normal subgroup of $\Gamma$. Let
$C_2=\GEN{\xi,\alpha_1,\dots,\alpha_n}$. By the construction of $\xi$, it follows that $C_2/C_1$ is central in
$\Gamma/C_1$. Since $t(\sigma,\tau) \in \GEN{\alpha_1,\dots,\alpha_n}\subseteq C_2$ for each $\sigma,\tau\in G$, one
has $\Gamma/C_2\simeq G$ and so it is enough to show that $G$ is supersolvable. Since $L/E_0$ is a Kummer extension,
$E_0(\alpha_i)/E_0$ is cyclic for every $i$. Moreover by the definition of $E_0$, $E_i=E_0(\alpha_1,\dots,\alpha_i)/K$
is Galois and hence $H_i=\Gal(L/E_i)$ is normal in $G$. Consider the composition series
    $$G \unrhd H_0\unrhd H_1 \unrhd H_2 \unrhd \cdots \unrhd H_{n-1} \unrhd H_n=1.$$
Since $G/H_0$ is abelian we only need to show that $H_i/H_{i+1}$ is cyclic for $0\le i\le n-1$. But $H_i/H_{i+1} \simeq
\Gal(E_{i+1}/E_i)=\Gal(E_i(\alpha_{i+1})/E_i)$ is an epimorphic image of $\Gal(E_0(\alpha_{i+1})/E_0)$ which is cyclic.
$\qed$\bigskip

To prove Theorem~\ref{RadSS}(b) we need the following lemma (see also \cite{AS1})   .

\begin{lemma}\label{NorSub}
Let $A=k(\Gamma)$ be a central simple $K$-algebra which is
projective Schur over $k$ and let $\Delta$ be a normal subgroup of
$\Gamma$. Then $B=k(\Delta)$ is semisimple. Moreover, if $B$ has
$n$ simple components then $A$ is isomorphic to $M_n(C)$, for some
algebra $C$ which is projective Schur over $k$. In fact $C$ is
spanned over $k$ by an epimorphic image of $\Gamma$.
\end{lemma}

\begin{proof}
Let $G=\Gamma k^\times/k^{\times}$ and let $\alpha\in H^2(G,k^{\times})$ be the cohomology class associated to the central extension
$1\rightarrow k^{\times} \rightarrow \Gamma k^\times \rightarrow G\rightarrow 1$. Then $A\simeq k^{\alpha} G/I$ for some
maximal ideal $I$ of $k^{\alpha} G$. Furthermore, putting $N=\Delta k^{\times}/k^{\times}$ and $\beta$ the restriction of $\alpha$ to
$N$, we have $B\simeq (k^{\beta} N \cdot I)/I \simeq k^{\beta} N/(k^{\beta} N\cap I)$. Thus $J(k^{\beta} N)\subseteq
J(k^{\alpha} G) \subseteq I$ (see \cite[Theorem 4.2]{Passman} or \cite[Corollary 2.5.30]{Rowen}) and hence $B$ is
semisimple.

Let $e$ be a primitive central idempotent of $B$ and $\{\sigma_1,\dots,\sigma_m\}$ a right transversal of the
centralizer of $e$ in $\Gamma$. Set $e_i=\sigma_i\inv g \sigma_i$ for each $i$. Since $\Delta$ is normal in $\Gamma$,
the elements $e_1,\dots,e_m$ are different primitive central idempotents of $B$ and hence $\sum_{i=1}^m e_i$ is a
non-zero central idempotent of $A$. Thus $\sum_{i=1}^m e_i=1$ since $A$ is simple. It follows that $B$ has $m$ simple
components and so $m=n$. Now, for every $g\in \Gamma$, the map $Ae \rightarrow Ag\inv e g = Aeg$ given by $x\mapsto xg$
is an isomorphism of left $A$-modules. Consequently, $A \simeq (Ae)^n$ as left $A$-modules and hence $A\simeq \End_A(A)
= M_n(\End(Ae)) = M_n(eAe)$. Finally, $eAe=k(e\Gamma e)$ is projective Schur over $k$ and the result follows.
\end{proof}

\noindent{{\bf Proof of Theorem~\ref{RadSS}(b)}. Let $A=k(\Gamma)$ be projective Schur over $k$ with $\Gamma$
supersolvable. By Lemma~\ref{NorSub}, we may assume that $k(\Delta)$ is simple for every normal subgroup $\Delta$ of
$\Gamma$. In particular, $k(\Delta)$ is a field if $\Delta$ is abelian.

We claim that $\Gamma'$, the commutator subgroup of $\Gamma$, is cyclic. Since $\Gamma$ is center by finite, by a
Theorem of Schur, $\Gamma'$ is finite (see e.g. \cite{Brown} Chapter IV, section 3, exercise 4(b)), and so there is a
subgroup $C$ of $\Gamma$ which is maximal among the cyclic normal subgroups of $\Gamma$ contained in $\Gamma'$. Let us
show that $\Gamma'=C$. If not, since $\Gamma$ is supersolvable, there is $a\in \Gamma'\setminus C$ such that
$D=\GEN{C,a}$ is normal in $\Gamma$. Furthermore, since $\Aut(C)$ is abelian, $\Gamma'$ centralizes $C$ and hence $D$
is abelian. It follows that $k(D)$ is a field and therefore $D$ must be cyclic. This contradicts the maximality of $C$
and the claim is proved.

Let $\Delta$ be a maximal abelian subgroup of $\Gamma$ which contains $\Gamma'(k^{\times}\cap \Gamma)$ (it exists since
$[\Gamma:(k^{\times}\cap \Gamma)]<\infty$). Clearly $\Delta$ is normal in $\Gamma$ and $G=\Gamma/\Delta$ is finite abelian.
Moreover, by the maximality of $\Delta$, $G$ acts faithfully on the field $L=k(\Delta)$. We obtain that $A\simeq
(L/K,\alpha)$ where $K$ is the field of $G$-invariant elements of $L$ and $\alpha$ is the cohomology class associated
to the group extension $1\rightarrow \Delta \rightarrow \Gamma \rightarrow G \rightarrow 1$. It is clear that $L/k$ is
a radical field extension and $\alpha\in Rad_k(L)$. Moreover $G\simeq \Gal(L/K)$ is abelian and therefore $A$ is
radical abelian over $k$ as desired.$\qed$

\section{The general case}\label{ProRep}

In this section we complete the proof of Theorem~\ref{Main} by proving (1 $\Rightarrow$ 3). Note that
Theorem~\ref{RadSS}(b) implies (1 $\Rightarrow$ 3) for projective Schur algebras $k(\Gamma)$ with $\Gamma$
supersolvable. In Proposition~\ref{RedE} below we reduce the general case to the supersolvable case however not over
the "right" field. Then, a suitable application of the corestriction map gives us the desired result. It is here, in
this very last step of the proof where we need the notion introduced above, namely of a central simple $K$-algebra
which is radical over $k$.

Recall that a semisimple $K$-algebra is split if it is isomorphic to a direct product of matrix algebras over $K$. If
$F/K$ is a field extension and $A$ is a semisimple $K$-algebra then we say that $F$ splits $A$ (or $F$ is a splitting
field of $A$) if $F\otimes_K A$ is a split $F$-algebra. Next lemma is basically known.

\begin{lemma}\label{SF}
Given a twisted group algebra $k^t G$ there exists a Galois
radical extension $L/k$ which splits $k^t H$ for every subgroup
$H$ of $G$. Consequently, if $H$ is any subgroup of $G$, then the
centre of any simple component of $k^t H$ is contained in $L$.
\end{lemma}

\begin{proof}
Let $n$ be the order of the element of $H^2(G,k^{\times})$ represented by the $2$-cocycle $t$. Then there is a $1$-coboundary
$\lambda:G\rightarrow k^{\times}$ such that $t(g,h)^n = \lambda_{g}\lambda_{h}\lambda_{gh}\inv$, for every $g,h\in G$. Let $F$
be the normal closure of $k(\lambda_{g}^{1/n}:g\in G)$ over $k$. Then $t$ is equivalent in $H^2(G,F^{\times})$ to the
$2$-cocycle $t'(g,h)=t(g,h)\lambda_{g}^{-1/n}\lambda_{h}^{-1/n}\lambda_{gh}^{1/n}$. Clearly $t'(g,h)^n=1$ and thus
$t'(g,h)$ is a root of unity, for every $g,h\in G$. Let $Z$ be the subgroup of $F^{\times}$ generated by the values
$\{t'(g,h):g,h\in G\}$ and let $G_1$ be the group extension $t':1\rightarrow Z \rightarrow G_1 \rightarrow G
\rightarrow 1$. It is easy to see that $F^t G$ is an epimorphic image of the group algebra $F G_1$ and hence, by the
Brauer splitting theorem there is a cyclotomic field extension $L/F$ which splits $k^t H$ for every subgroup $H$ of
$G$.

If $A$ is a simple epimorphic image of $k^t H$ then $A\otimes_k L$ is an epimorphic image of $k^t H\otimes_k L$ and
therefore $A\otimes_k L$ is a direct product of matrix algebras over $L$. Then $Z(A)\otimes_k L = Z(A\otimes_k L) =
L^m$ for $m=[Z(A):k]$ and since $L/k$ is a Galois extension, $Z(A)\subseteq L$. This proves the lemma.
\end{proof}

\begin{proposition}\label{RedE}
Let $k^t G$ be a twisted group algebra over a field $k$ of characteristic zero. Let $L$ be a finite radical extension
of $k$ which splits $k^t H$ for every finite subgroup $H$ of $G$. Let $A$ be an epimorphic image of $k^t G$ with centre
$K$. Assume that the index of $A$ is a power of a prime integer $p$. Then there is a supersolvable subgroup $H$ of $G$,
a $K'$-central simple algebra $B$ and a subfield $E$ of $L$ such that:

a) $B$ is an epimorphic image of $k^t H$.

b) $E$ contains $K$ and $K'$.

c) $[E:K]$ is coprime with $p$.

d) $E\otimes_K A$ is equivalent to the $p$-th part of
$E\otimes_{K'} B$ in $Br(E)$.

\end{proposition}

The proof of Proposition~\ref{RedE} requires the theory of
projective characters. For the reader's convenience we recall some
of the tools used in the proof. These tools can be founded in
\cite{CR}, \cite{Isaacs} and \cite{Kar}.

Given a twisted group algebra $k^t G$ we fix a $k$-basis $\{u_g:g\in G\}$ where $u_g u_g = t(g,h) u_{gh}$ for every
$g,h\in G$. We also fix a radical extension $L$ of $k$ which splits $k^t H$ for every subgroup $H$ of $G$
(Lemma~\ref{SF}) and consider $\{u_g:g\in G\}$ as an $L$-basis of $L^t G$.

If $M$ is a left $L^t G$-module then the projective
$t$-representation associated to $M$ is the map $\rho:G
\rightarrow \GL_L(M)$ that associates to $g\in G$ the action of
$u_g$ on $M$. The projective $t$-representation $\rho$ satisfies
the following relation:
    \begin{equation}\label{PRC}
    \rho(g)\rho(h) = t(g,h) \rho(gh).
    \end{equation}
Conversely, if $M$ is an $L$-vector space, then every map
$\rho:G\rightarrow \GL_L(M)$ satisfying (\ref{PRC}) gives an $L^t
G$-module structure on $M$ whose projective $t$-representation is
$\rho$. The composition of $\rho$ with the trace map
$\tr:\GL_L(M)\rightarrow L$ is called the projective $t$-character
afforded by $M$ or by $\rho$.  Two left $L^t G$-modules are
isomorphic if and only if they afford the same projective
$t$-character \cite[7.1.11]{Kar}. A projective $t$-representation
(resp. projective $t$-character) is said to be irreducible if it
is the projective $t$-representation of (resp. the projective
$t$-character afforded by) an irreducible $L^t G$-module.


Let $F$ be a subfield of $L$ containing $k$. Every irreducible projective $t$-character $\chi$ induces an $L$-linear
map $\overline{\chi}:L^t G\rightarrow L$ given by $\overline{\chi}(u_g)=\chi(g)$ and there is a unique simple component
$A$ of $F^t G$ such that $\overline{\chi}(A)\ne 0$. Conversely for every simple component $A$ of $F^t G$ there is a
(non-necessarily unique) irreducible projective $t$-character $\chi$ of $G$ such that $\overline{\chi}(A)\ne 0$.
Furthermore, if $F=Z(A)$ then the character $\chi$ such that $\overline{\chi}(A)\ne 0$ is unique and it is given by
$\chi(g)=\tr(\pi(u_g))$, where $\pi:F^t G\rightarrow A$ is the projection and $\tr:A\rightarrow F$ is the reduce trace
of $A$. We will refer to this $\chi$ as the projective $t$-character given by the reduced trace of $A$.

Given $F$ as above and $\chi$ a projective $t$-character, $F(\chi)$ denotes the extension of $F$ generated by
$\{\chi(g):g\in G\}$. One says that $\chi$ is representable over $F$ if $\chi$ is the projective $t$-representation
associated to $M_L=L\otimes_F M$ for some left $F^t G$-module $M$. Equivalently $F$ splits the unique simple component
of $k^t G$ on which $\overline{\chi}$ does not vanish. Clearly if a projective $t$-character is representable over $F$
then it takes values in $F$ but the converse is not true in general.

Assume now that $\chi$ is irreducible and let $A$ be the unique simple component of $F^t G$ such that
$\overline{\chi}(A)\ne 0$. The Schur index $m_F(\chi)$ of $\chi$ over $F$ is by definition $\Ind(A)$, the index of $A$.
Equivalently, it may be defined by $m_F(\chi) = \min \{[E:F(\chi)]: \chi \mbox{ is representable over } E \}$ (see
\cite[8.3.1]{Kar} or \cite[13.4]{Pie}). Furthermore the field $F(\chi)$ is $F$-isomorphic (but not necessarily equal)
to $Z(A)$ \cite[8.2.4]{Kar} and hence there is $\sigma\in \Gal(L/F)$ such that $F(\sigma\circ \chi)=Z(A)$.

Given two projective $t$-characters $\chi_1$ and $\chi_2$ afforded
by left $L^t G$-modules $M_1$ and $M_2$ one defines the scalar
product of $\chi_1$ and $\chi_2$ as
    $$(\chi_1,\chi_2)=\dim_L \Hom_{L^t G}(M_1,M_2)$$
Let us recall some properties of this product. Let $\chi$ and
$\lambda$ be projective $t$-characters afforded by modules $M$ and
$N$ and assume that $\chi$ is irreducible:
\begin{enumerate}
\item $(\chi,\lambda)$ is the number of direct summands isomorphic
to $M$ in a (any) decomposition of $N$ as direct
    sum of simple $L^t G$-modules.

    In particular, if $\lambda$ is irreducible then $(\chi,\lambda)=1$
if $\chi=\lambda$ and $0$ otherwise. \item The irreducible
projective $t$-characters are linearly independent over $L$. \item
If $\lambda$ is representable over $F$ then $m_F(\chi)$ divides
$(\chi,\lambda)$.
\end{enumerate}

For any subgroup $H$ of $G$ we denote by $t_H$ and by $\chi_H$ the
restrictions of $t$ and  $\chi$ to $H$. Note that $\chi_H$ is a
projective $t_H$-character. On the other hand, if $\theta$ is the
projective $t_H$-character afforded by the left $L^{t_H} H$-module
$M$, then the projective $t$-character $\theta^G$ induced from
$\theta$ to $G$ is the character afforded by $L^t G \otimes_{L^t
H} M$ as left $L^t G$-module.

Let $k^{t_1} G_1$ and $k^{t_2} G_2$ be two twisted group algebras.
Let $L/k$ be a field extension that splits $k^{t_1} G_1$ and
$k^{t_2} G_2$. Let $M_i$ be left $L^{t_i}G_i$-modules for $i=1,2$
and let $\chi_i$ be the corresponding characters. Then the map
    $$(t_1 \times t_2)((g_1,g_2),(h_1,h_2)) = t_1(g_1,h_1)t_2(g_2,h_2)$$
is a $2$-cocycle on $G_1\times G_2$. Furthermore, since $L^{t_1} G_1 \otimes_L L^{t_2} G_2 \simeq L^{t_1\times t_2}
(G_1\times G_2)$ we have that $M_1\otimes_L M_2$ is a left $L^{t_1\times t_2}(G_1 \times G_2)$-module and we denote its
character by $\chi_1\otimes \chi_2$. In the particular case where $G=G_1=G_2$ we may identify $G$ with $\{(g,g):g\in
G\}$ and view $M_1\otimes_L M_2$ as a left $L^{(t_1\times t_2)_{G}} G$-module. Note that $(t_1\times t_2)_{G}$ is just
the pointwise multiplication $t_1t_2$. We denote the character $(\chi_1\otimes \chi_2)_{G}$ by $\chi_1\cdot \chi_2$.

For the proof of Proposition~\ref{RedE} we need the following
properties.

\begin{enumerate}
\setcounter{enumi}{3}
\item
Let $k^{t_1} G_1$, $k^{t_2} G_2$, $\chi_1$ and $\chi_2$ be as above.
\begin{enumerate}
\item
$(\chi_1\otimes \chi_2)(g_1,g_2)=\chi_1(g_1)\chi_2(g_2)$ for each $(g_1,g_2)\in G_1\times G_2$. \item If $\lambda_i$ is
a projective $t_i$-character ($i=1,2$) then $(\chi_1\otimes \chi_2,\lambda_1\otimes \lambda_2) =
(\chi_1,\lambda_1)(\chi_2,\lambda_2)$. In particular, $\chi_1\otimes \chi_2$ is irreducible if and only if $\chi_1$ and
$\chi_2$ are irreducible.
\item
Assume that $\chi_1$ and $\chi_2$ are irreducible and for $i=1,2$ let $A_i$ be a simple component of $k^{t_i} G_i$ with
centre $F$ such that $\overline{\chi_i}(A_i)\ne 0$. Then $\overline{\chi_i\otimes \chi_2}(A_1\otimes_F A_2)\ne 0$.
\end{enumerate}
\item Let $H$ be a subgroup of $G$, $\chi$ a projective
$t$-character and $\theta$ a projective $t_H$-character.
\begin{enumerate}
\item $\chi \cdot \theta^G = (\chi_H \cdot \theta)^G$
\cite[5.5.1]{Kar}. \item (Frobenius Reciprocity) If $\chi$ and
$\theta$ are irreducible then $(\chi_H,\theta)=(\chi,\theta^G)$
    \cite[5.6.3]{Kar}.
\end{enumerate}
\end{enumerate}


\noindent{\bf Proof of Proposition~\ref{RedE}} We will need the following weak version of the Berman-Witt Theorem
\cite[Section 2.21A]{CR} for ordinary characters: There exist supersolvable subgroups $H_1,\dots,H_n$ of $G$,
characters $\theta_1,\dots,\theta_n$ of the $H_i$'s, with $\theta_i(H_i)\subseteq K$ and integers $a_1,\dots,a_n$ such
that $1=\sum_i a_i \theta_i^G$.

The group $\Gal(L/K)$ acts on the set of irreducible $t_{H_i}$-characters by $\sigma\cdot \xi = \sigma\circ \xi$,
($\sigma\in \Gal(L/K)$, $\xi$ an irreducible $t_{H_i}$ character). Let $R_i$ be a set of representatives for the
corresponding orbits and for each $\xi\in R_i$ let $\tr(\xi)$ be the sum of the elements in the orbit of $\xi$. Let
$\chi$ be an irreducible projective $t$-character of $G$ such that $\chi(G)\subseteq K$ and $\overline{\chi}(A)\ne 0$
(e.g. take $\chi$ the projective $t$-character given by the reduced trace of $A$). Note that the projective
$t_{H_i}$-character $\chi_{H_i}\cdot \theta_i$ takes values in $K$ and hence $(\chi_{H_i}\cdot
\theta_i,\xi)=(\chi_{H_i}\cdot \theta_i,\sigma\cdot \xi)$ for every $\sigma\in \Gal(L/K)$. It follows that
$\chi_{H_i}\theta_i = \sum_{\xi\in R_i} b_{\xi} \tr(\xi)$ where the $b_{\xi}$'s are non-negative integers. Then we have

    $$
    \chi = \chi \cdot 1 = \sum_i a_i (\chi \cdot \theta_i^G) \\
    = \sum_i a_i (\chi_{H_i} \cdot \theta_i)^G \\
    = \sum_i \sum_{\xi\in R_i} a_i b_{\xi} \tr(\xi)^G$$
and thus
    $$1 = (\chi,\chi) = \sum_i \sum_{\xi\in R_i} a_i b_{\xi}
(\chi,\tr(\xi)^G) =
    \sum_i \sum_{\xi\in R_i} a_i b_{\xi}  [K(\xi):K](\chi,\xi^G).$$
It follows that there is an $i$ and a $\xi\in R_i$ such that $s=(\chi_{H_i},\xi)=(\chi,\xi^G)$ and $[E=K(\xi):K]$ are
both coprime with $p$, the unique (possible) prime divisor of the Schur index of $\chi$ over $K$ ($=\Ind(A)$). Put
$H=H_i$ and let $B$ be the unique simple component of $k^{t_{H}} H$ such that $\xi(B)\ne 0$. If necessary, by replacing
$\xi$ with another element in its orbit we may assume that $k(\xi)=Z(B)$ and therefore $Z(B)$ is contained in $E$. We
are to show that the algebra $B$ is equivalent to $A$ over the field $E$. To this end consider the $k$-linear map
$\phi:k^{t\inv} G \rightarrow k^t G$ given by $u_g \mapsto u_g\inv$. Note that $\phi$ is an anti-isomorphism of
$k$-algebras and so if $\pi:k^t G \rightarrow A$ is an epimorphism of $k$-algebras, $A^\circ$ is the opposite algebra
of $A$ and $1^\circ$ is the identity map on the underlying sets, then the composition $1^\circ \circ \pi \circ
\phi:k^{t\inv} G {\rightarrow} A^\circ$ is an epimorphism of $k$-algebras.
Let $\chi\inv$ be the projective $t\inv$-character given by the reduce trace of $A^\circ$. Then $\chi\otimes \chi\inv$
does not vanish on $A\otimes_K A^\circ$ which is split and hence $\chi\otimes \chi\inv$ is representable over $K$. This
implies that $\chi_H\otimes \chi\inv$ is representable over $K$ and so the Schur index of $\xi\otimes \chi\inv$ over
$E$ divides $(\chi_H\otimes \chi\inv,\xi\otimes \chi\inv)=(\chi_H,\xi)(\chi\inv,\chi\inv)=(\chi_H,\xi)=s$. Moreover
$C=B_E \otimes_E A^\circ_E$ is the unique simple component of $E^{t\times t\inv}(H\times G)$ such that $(\xi\otimes
\chi\inv)(C)\ne 0$. Since $s$ is coprime with $p$, the $p$-th part of $[C]$ is trivial. Finally, since $[B_E\otimes_E
A^\circ_E]=1=[A_E\otimes_E A^\circ_E]$, $A_E$ is equivalent to the $p$-th part of $B_E$ and we are done.
    $\qed$\bigskip


To complete the proof of (1 $\Rightarrow$ 3) (and henceforth of
Theorem~\ref{Main}) we need an additional lemma.

\begin{lemma}\label{TPR}
Let $A_1$ and $A_2$ be central simple $K$-algebras and $k$ a
subfield of $K$.
\begin{itemize}
\item[(a)] If $A_1$ and $A_2$ are projective Schur over $k$ then
so is $A_1\otimes_K A_2$. \item[(b)] If $A_1$ and $A_2$ are
radical (resp. radical abelian) over $k$ then $A_1\otimes_K A_2$
is equivalent in
    $Br(K)$ to an algebra which is radical (resp. radical
abelian) over $k$. \item[(c)] If $A=(F/K,\alpha=[t])$ is radical
over $k$ and $F_1/K$ is a Galois extension where $F_1$ is radical
over $k$ and contains $F$, then the inflation of $A$ to $F_1$ is a
radical over $k$. \item[(d)] If $A=A_1$ is radical over $k$ and
$E$ is a extension of $K$ contained in a radical extension of $k$
then
    $A_E$ is equivalent in $Br(E)$ to an algebra which is radical over $k$.
\end{itemize}
\end{lemma}

\begin{proof}
(a) is obvious.

(b) Let $A_i=(L_i/K,\alpha_i)$ be radical (resp. radical abelian) over $k$ for $i=1,2$. Then $L=L_1L_2$ is a radical
extension of $k$ (resp. radical extension of $k$ and $L/K$ is abelian). It follows that the inflated algebra
$B_i=(L/K,\bar{\alpha_i})$ is radical (resp. radical abelian) over $k$. Clearly $A_i$ is equivalent to $B_i$,
$B_1\otimes_K B_2$ is equivalent to $A=(L/K,\bar{\alpha_1}\bar{\alpha_2})$ and the latter is radical (resp. radical
abelian) over $k$.

(c) The inflation of $A$ to $F_1$ is
$\widehat{A}=(F_1/K,[\widehat{t}])$ where
$\widehat{t}(g_1,g_2)=t(\overline{g_1},\overline{g_2})$ and
$g\mapsto \overline{g}$ is the restriction homomorphism
$\Gal(F_1/K)\rightarrow \Gal(F/K)$ and so $\widehat{A}$ is radical
over $k$.

(d) Assume that $A=(F/K,\alpha)$ is radical over $k$ and let $F_1$
be a radical extension of $k$ containing $E$ and $F$. Then $A_E$
is equivalent to $B_E$ where $B$ is the inflation of $A$ to $F_1$.
By 3, $B$ is radical over $k$, and therefore its restriction $B_E$
is radical over $k$.
\end{proof}

\noindent{{\bf Proof of (1 $\Rightarrow$ 3) of Theorem~\ref{Main}}. Let $A$ be a central simple $K$-algebra which is
projective Schur over $k$ and let $k^{t} G {\rightarrow} A$ be an epimorphism of $k$-algebras. We show that $A$ is
radical over $k$.

{\em Claim}: We may assume that the index of $A$ is a prime power.

Indeed, if $p_1,\dots,p_n$ are the prime divisors of the index of $A$ (equivalently, of the order of the class of $A$
in $Br(K)$), then $A$ is equivalent in $Br(K)$ to an algebra of the form $A_1\otimes_K \cdots \otimes_K A_n$ where
$A_i$ has index (and order) a power of $p_i$ for each $i$. Then for each $i$ there is a positive integer $m$ such that
$A_i$ is equivalent to $A^{\otimes m}$ and hence by Lemma~\ref{TPR}(a) $A_i$ is projective Schur over $k$. If (1
$\Rightarrow$ 3) holds for algebras of prime power index then $A_i$ is equivalent to an algebra $B_i$ which is radical
abelian over $k$. Applying Lemma~\ref{TPR}(b) one obtains that $B_1\otimes \cdots \otimes B_n$ equivalent to an algebra
which is radical abelian over $k$ and the claim follows.

So let $\Ind(A)$ be a power of a prime $p$. By Proposition~\ref{RedE} there exist a supersolvable subgroup $H$ of $G$,
a simple epimorphic image $B$ of $k^t H$ and a field $E$ containing $K=Z(A)$ and $Z(B)$, such that $[E:K]$ is not a
multiple of $p$ and $A_E$ is equivalent to the $p$-th part of $B_E$. Applying Theorem~\ref{RadSS}(b), $B$ is radical
abelian over $k$ and hence by Lemma~\ref{TPR}(d) $B_E$ is radical (non-necessarily abelian) over $k$. Clearly, $A_E$ is
equivalent (in $Br(E)$) to $B_E^{\otimes m}$ for some integer $m$ and hence, by Lemma~\ref{TPR}(b)}, is equivalent to
an algebra $C=(F/E,\alpha)$ which is radical $k$. Inflating $C$ to the normal closure of $F$ over $K$ we may assume
that $F/K$ is Galois (and radical over $k$).

Let $Q=\Gal(F/K)$, $T=\Gal(F/E)$ and let $cor: H^2(T,F^{\times}) \rightarrow H^2(Q,F^{\times})$  be the corestriction map. Clearly,
the Galois group $Q$ acts on $Rad_k(F)$. Since $C$ is radical over $k$, $\alpha\in Rad_k(F)$ (i.e. $\alpha$ belongs to
the image of $H^2(T,Rad_k(F))\rightarrow H^2(T,F^{\times})$) and hence $cor(\alpha)\in Rad_k(F)$ by the commutativity of the
diagram
    $$\matriz{{cccccc}
    H^{2}(T,Rad_k(F))      & \rightarrow  & H^{2}(T, F^{\times}) \\
    cor \downarrow  &              & cor \downarrow \\
    H^{2}(Q, Rad_k(F))     &  \rightarrow & H^{2}(Q, F^{\times})}$$
It follows that $cor(C)=(F/K,cor(\alpha))$ is a central simple $K$-algebra radical over $k$. Now, since $C$ is
equivalent to $A_E$, we have that $cor(C)$ is equivalent to $A^{\otimes n}$ where $n=[Q:T]=[E:K]$. Furthermore, since
$n$ is coprime with $p$, $A$ is equivalent to $cor(C)^{\otimes m}$ for some integer $m$. By Lemma~\ref{TPR}(b),
$cor(C)^{\otimes m}$ is equivalent to a radical algebra over $k$. This proves (1 $\Rightarrow$ 3) and completes the
proof of Theorem~\ref{Main}.$\qed$

\noindent Department of Mathematics, Technion-Israel Institute of
Technology, 32000 Haifa, Israel,
\\ aljadeff@techunix.technion.ac.il

\noindent Departamento de Matem\'{a}ticas, Universidad de Murcia, 30100 Murcia, Spain, adelrio@um.es

\end{document}